\documentclass[12pt]{amsart}
\usepackage{amsmath,amscd,amssymb,amsfonts,graphics,mathrsfs}
\usepackage{hyperref}
\newcommand{\hl}{\hyperlink}
\newcommand{\htt}{\hypertarget}
\newcommand{\h}{\hbox}
\newcommand{\q}{\quad}
\newcommand{\nin}{\noindent}
\newcommand{\bs}{\par\bigskip}
\newcommand{\ms}{\par\medskip}
\newcommand{\sk}{\par\smallskip}
\newcommand{\bsn}{\par\bigskip\noindent}
\newcommand{\msn}{\par\medskip\noindent}
\newcommand{\skn}{\par\smallskip\noindent}
\newcommand{\ges}{\geqslant}
\newcommand{\les}{\leqslant}
\newcommand{\1}{\hskip1pt}
\newcommand{\mcap}{\hbox{$\bigcap$}}

\newcommand{\mopl}{\hbox{$\bigoplus$}}
\newcommand{\msum}{\hbox{$\sum$}}

\newcommand{\mwedge}{\hbox{$\bigwedge$}}
\newcommand{\A}{{\mathscr A}}

\newcommand{\D}{{\mathscr D}}

\newcommand{\Hc}{{\mathscr H}}
\newcommand{\Lc}{{\mathscr L}}
\newcommand{\Nc}{{\mathscr N}}
\newcommand{\I}{{\mathscr I}}
\newcommand{\OO}{{\mathscr O}}
\newcommand{\M}{{\mathscr M}}
\newcommand{\Hom}{\mathscr{H}\!om}
\newcommand{\Mct}{\widetilde{\mathscr M}}
\newcommand{\Mt}{\widetilde{M}}
\newcommand{\Xt}{\widetilde{X}}
\newcommand{\Yt}{\widetilde{Y}}
\newcommand{\xt}{\widetilde{x}}
\newcommand{\tti}{\widetilde{t}}

\newcommand{\PP}{{\mathbb P}}
\newcommand{\Q}{{\mathbb Q}}
\newcommand{\C}{{\mathbb C}}
\newcommand{\Ab}{{\mathbb A}}

\newcommand{\Z}{{\mathbb Z}}
\newcommand{\Sp}{{\rm Sp}}
\newcommand{\Gr}{{\rm Gr}}

\newcommand{\al}{\alpha}
\newcommand{\be}{\beta}

\newcommand{\la}{\lambda}
\newcommand{\te}{\theta}
\newcommand{\om}{\omega}
\newcommand{\Om}{\Omega}

\newcommand{\dd}{\partial}
\newcommand{\ddd}{{\rm d}}
\newcommand{\eq}{\,{=}\,}
\newcommand{\co}{\,{:}\,}
\newcommand{\gess}{\,{\ges}\,}
\newcommand{\less}{\,{\les}\,}
\newcommand{\sgt}{\,{>}\,}
\newcommand{\slt}{\,{<}\,}
\newcommand{\nes}{\,{\ne}\,}
\newcommand{\notins}{\,{\notin}\,}
\newcommand{\mi}{\1{-}\1}

\newcommand{\pl}{\1{+}\1}

\newcommand{\bl}{\bigl}
\newcommand{\br}{\bigr}
\newcommand{\sst}{\,{\subset}\,}
\newcommand{\stm}{\,{\setminus}\,}

\newcommand{\ins}{\,{\in}\,}
\newcommand{\tos}{\,{\to}\,}
\newcommand{\defs}{\,{:=}\,}
\newcommand{\ssc}{\,\raise.15ex\hbox{${\scriptstyle\circ}$}\,}
\newcommand{\ssb}{\raise.15ex\h{${\scriptscriptstyle\bullet}$}}
\newcommand{\into}{\hookrightarrow}

\newcommand{\simto}{\,\,\rlap{\hskip1.3mm\raise1.4mm\hbox{$\sim$}}\hbox{$\longrightarrow$}\,\,}
\begin{document}
\h{}\bs\bs
\centerline{\large Verdier Specialization and Restrictions of Hodge modules}
\bs
\centerline{Qianyu Chen, Bradley Dirks, Morihiko Saito}
\bs\bs
\vbox{\nin\narrower\smaller
{\bf Abstract.} We give an explicit formula to express the cohomological pullback functors of Hodge modules under closed immersions of smooth varieties using Verdier specializations and $V$-filtrations of Kashiwara and Malgrange. This was locally obtained by the first two authors assuming the existence of global defining functions. We also give a quite simplified proof of the theorem reducing to the monodromical case via the Verdier specialization and using induction on codimension.}
\bs\bs\ms
\centerline{\bf Introduction}
\bsn
For a closed immersion of smooth complex varieties $i_Z\co Z\into X$, the cohomological pullback functors $H^{\ssb}i_Z^*\M$, $H^{\ssb}i_Z^!\M$ for a mixed Hodge module $\M$ on $X$ can be calculated locally by iterating the nearby and vanishing cycle functors along local coordinates $x_1,\dots,x_r$ such that $Y\eq\mcap_{i=1}^r\{x_i\eq0\}$, see for instance \cite[Remark after Cor.\,2.24]{mhm}. Here we can replace $\M$ with its {\it Verdier specialization\1} $\Sp_Z\M$ and $i_Z$ with the immersion to the normal bundle of $Z\sst X$ as the zero-section, see \cite[2.30]{mhm} (locally), Proposition~\hl{P1.2}{1.2} and also (\hl{7}{7}) below. We may thus assume that $X$ is a vector bundle over $Z$ and the Hodge module is monodromical. (In this paper we say that a mixed Hodge module on $X$ is {\it monodromical\1} if the isomorphism (\hl{2}{2}) below holds for the underlying $\D$-module in a compatible way with the Hodge and weight filtrations. The last condition is satisfied for the Verdier specialization of a mixed Hodge module, see \hl{1.2}{1.2} below.)
\sk
Let $\te$ be the vector field on $X$ corresponding to the $\C^*$-action on the vector bundle. Its action on $\I_Z/\I_Z^2$ is the identity, where $\I_Z\sst\OO_X$ is the ideal sheaf of $Z\sst X$. Let $(M';F,W)$ be the underlying bifiltered left $\D_X$-module of a monodromical mixed Hodge module $\M'$, for instance the Verdier specialization $\Sp_Z\M$ of a mixed Hodge module $\M$. Set
\htt{1}{}
$$M'_{\al}\defs{\rm Ker}\bl((\te\mi\al\pl r)^k:M'|_Z\tos M'|_Z\br)\q(k\gg0,\,\al\ins\Q).
\leqno(1)$$
(In this paper $|_Z$ denotes the sheaf-theoretic restriction.)
\sk
We now assume that the vector bundle $\pi\co X\tos Z$ is trivial (shrinking $Z$ if necessary). Choosing a trivialization, we have the isomorphism as $\D_Z$-modules:
\htt{2}{}
$$\pi_*M'\eq\mopl_{\al\in\Q}\,M'_{\al},
\leqno(2)$$
since $\M'$ is monodromical. The $M'_{\al}$ are {\it not\1} coherent $\D_Z$-modules unless ${\rm Supp}\,M'\sst Z$ or $r\eq1$.
\sk
Let $t_1,\dots,t_r$ be sections of the dual vector bundle trivializing the bundle. These are identified with functions on $X$ which are linear on each fiber. Let $\dd_i$ be the vector fields defined on $X$ such that their integral curves are contained in the fibers of $\pi$ and $\langle\dd_i,\ddd\1t_j\rangle\eq\delta_{i,j}$ ($i,j\ins[1,r]$). We have $\te\eq\msum_{i=1}^r\,t_i\dd_i$ on $X$.
\sk
Consider the Koszul complexes $K_{\!Z}^{!,\ssb}(M')$ and $K_{\!Z}^{*,\ssb}(M')$ for the actions of $t_1,\dots,t_r$ and $\dd_1,\dots,\dd_r$ on $\pi_*M'$ respectively. They have the Hodge filtration $F$ together with the $\Q$-grading such that $F_pK_{\!Z}^{!,j}(M')_{\al}$ and $F_pK_{\!Z}^{*,j}(M')_{\al}$ are respectively direct sums of
\htt{3}{}
$$F_{p+r}M'_{\al+j}\q\h{and}\q F_{p+j}M'_{\al+r-j}\q(\al\ins\Q,\,p\ins\Z,\,j\ins[0,r]).
\leqno(3)$$
These can be made independent of the choice of a local trivialization given by the $t_i$, using the exterior products of the conormal sheaf $\Nc_{Z/X}^*\defs\I_Z/\I_Z^2$ and considering
\htt{4}{}
$$\mwedge^{r-j}\Nc^*_{Z/X}{\otimes}_{\OO_Z}F_{p+r}M'_{\al+j}\q\h{and}\q\mwedge^j\Nc^*_{Z/X}{\otimes}_{\OO_Z}F_{p+j}M'_{\al+r-j},
\leqno(4)$$
instead of the above direct sums. Here the differentials are given respectively by
\htt{5}{}
$$\aligned\eta{\otimes}\xi\mapsto\msum_{p=1}^{r-j}\,\msum_{k=1}^r\,(-1)^{p-1}\eta^{(p)}{\otimes}\langle\dd_k,\eta_p\rangle t_k\xi,&\\\h{and}\,\q\q\q\q\q\q\eta{\otimes}\xi\mapsto\msum_{k=1}^r\,\ddd t_k{\wedge}\eta{\otimes}\dd_k\xi,&\endaligned
\leqno(5)$$
for $\xi\ins M'$, $\eta\eq\bigwedge_{i=1}^{j'}\,\eta_i$ with $\eta_i\ins\Nc^*_{Z/X}$ ($j'\eq r{-}j$ or $j$), and $\eta^{(p)}\defs\bigwedge_{i\ne p}\,\eta_i$. We denote by $\ddd t_k$ the image of $t_k$ in $\Nc^*_{Z/X}$. We get the initial Koszul complexes if $\eta_i\eq\ddd t_{k_i}$ for any $i$. We can show the well-definedness of the action of $\D_Z$ on the cohomology sheaves of the Koszul complexes using a {\it flat or injective resolution of $M'$ as a $\D_X$-module.} (Indeed, the ambiguities of liftings of vector fields on $Z$ are given by $\OO_Z$-linear combinations of the $t_i\dd_j$ when a trivialization of vector bundle is changed, see Remark~\hl{R2.1d}{2.1d} below.) For $K_{\!Z}^{!,\ssb}(M')$, this is closely related to \cite[(1.1.12)]{RSW}. The modified Koszul complexes are defined globally. Forgetting the Hodge filtration and the grading, these Koszul complexes are isomorphic to the pullbacks $i^!_ZM'$ and $i^*_ZM'[-r]$ as $\D$-modules (see Remark~\hl{R2.1e}{2.1e}--\hl{R2.1f}{f} below) where $i^*_ZM'$ can be identified with the direct image of $M'$ as $\D$-module by $\pi$ in the monodromical case as is well known. It is rather easy to show the acyclicity of $K_{\!Z}^{!,\ssb}(M')_{\al}$, $K_{\!Z}^{*,\ssb}(M')_{\al}$ for $\al\nes0$, see \cite{CD} and also (\hl{9}{9}) below. The following gives a globalization of an assertion in \cite{CD} (whose proof is quite complicated) in the monodromical Hodge module case.
\par\htt{T1}{}\msn
{\bf Theorem 1.} {\it The Hodge filtration $F$ on the degree $0$ part of the Koszul complexes $K_{\!Z}^{!,\ssb}(M')_0$, $K_{\!Z}^{*,\ssb}(M')_0[r]$ are strict. Their $F$-filtered cohomology $\D_Z$-modules have the weight filtration $W$ induced from the relative monodromy filtration $W'$ for the actions of $\te\mi\al\pl r$ on the $W$\!-filtered $\D_Z$-module $M'_{\al}$ using $K_{\!Z}^{!,\ssb}(W'_{\ssb}M')_0$, $K_{\!Z}^{*,\ssb}(W'_{\ssb}M')_0[r]$ and shifted by cohomology degree. These bifiltered $\D_Z$-modules are isomorphic to the bifiltered underlying $\D_Z$-modules of the cohomological pullbacks of mixed Hodge modules $H^{\ssb}(i'_Z)^!\M'$ and $H^{\ssb}(i'_Z)^*\M'$ respectively, where $i'_Z\co X\into T_ZX$ is the inclusion.}
\ms
Here $H^{\ssb}\co D^b\A\tos\A$ denotes the standard cohomological functor with $\A\eq{\rm MHM}(Z)$. In the monodromical case the relative monodromy filtration is the same as the original weight filtration, see Remark~\hl{R2.1a}{2.1a} below. Theorem~\hl{T1}{1} actually implies a globalization of the assertion in \cite{CD} mentioned before Theorem~\hl{T1}{1} as follows.
\par\htt{T2}{}\msn
{\bf Theorem 2.} {\it Let $X$ be a smooth complex variety, and $Z\sst X$ be a smooth closed subvariety of codimension $r$. Let $\M$ be a mixed Hodge module on $X$ with $(M;F,W)$ the underlying bifiltered $\D_X$-module. Let $V$ be the $V$-filtration of Kashiwara and Malgrange on $M$ along $Z$ indexed by $\Q$ $($instead of $\Z)$. Then we have the isomorphisms of filtered $\D_Z$-modules
\htt{6}{}
$$(M'_{\al},F)\eq\Gr_V^{\al}(M,F)\q(\forall\,\al\ins\Q),
\leqno(6)$$
with $(M',F)$ the underlying filtered $\D$-module of $\M'\defs\Sp_Z(\M)$, see {\rm(\hl{1.2.8}{1.2.8})} below. Defining the $F$-filtered Koszul complexes as in {\rm(\hl{3}{3})}, the assertion of Theorem~{\rm\hl{T1}{1}} holds with $H^{\ssb}(i'_Z)^!\M'$, $H^{\ssb}(i'_Z)^*\M'$ replaced by $H^{\ssb}i_Z^!\M$, $H^{\ssb}i_Z^*\M$ respectively with $i_Z\co Z\into X$ the inclusion. The weight filtration $W$ is induced from the relative monodromy filtration $W'$ for the action of $\theta{-}\al{+}r$ on $\Gr_V^{\al}(M,W)$ using $K_{\!Z}^{!,\ssb}(W'_{\ssb}M')_0$, $K_{\!Z}^{*,\ssb}(W'_{\ssb}M')_0[r]$, and is shifted by cohomology degree.}
\ms
Indeed, we have the isomorphisms in $D^b{\rm MHM}(Z)$
\htt{7}{}
$$i^*_Z\M=(i'_Z)^*\M',\q i^!_Z\M=(i'_Z)^!\M'\q\h{with}\,\,\,\M'\defs\Sp_Z\M,
\leqno(7)$$
see Proposition~\hl{P1.2}{1.2} below. Moreover the relative monodromy filtration coincides with the one associated with the Verdier specialization, see \hl{1.2}{1.2} below. So Theorem~\hl{T2}{2} follows from Theorem~\hl{T1}{1}. If we do {\it not\1} use the Verdier specialization, we would have to assume the existence of global defining functions $t_1,\dots,t_r$ of $Z\sst X$ together with {\it mutually commuting vector fields\1} $\dd_1,\dots,\dd_r$ on $X$ with $\langle\dd_i,\ddd t_j\rangle\eq\delta_{i,j}$ for any $i,j\ins[1,r]$, see \cite{CD}. Here the mutual commutativity does {\it not\1} follow from the last condition unless the $t_i$ are extended to local coordinates of $X$ and we take the associated vector fields. Note also that the existence of the $\dd_i$ does not follow from that of the $t_i$ even in the case $r\eq1$ unless we have a trivial family in the proper case. (We can, however, use the graph embedding by the $t_i$ to get the vector fields.)
\sk
By increasing induction on $r$, we prove Theorem~\hl{T1}{1} together with the following filtered and unfiltered acyclicities (which are also proved in \cite{CD}):
\htt{8}{}
$$\aligned\Hc^j\Gr^F_pK_{\!Z}^{!,\ssb}(M')_{\al}&\eq0\q\h{if}\,\,\,\al\sgt0,\\ \Hc^j\Gr^F_pK_{\!Z}^{*,\ssb}(M')_{\al}&\eq0\q\h{if}\,\,\,\al\slt0.\endaligned
\leqno(8)$$
\vskip-3mm
\htt{9}{}
$$\Hc^jK_{\!Z}^{!,\ssb}(M')_{\al}\eq\Hc^jK_{\!Z}^{*,\ssb}(M')_{\al}\eq0\q\h{if}\,\,\,\al\nes0.
\leqno(9)$$
\sk
Theorem~\hl{T1}{1} is closely related to Griffiths' theorem on rational integrals \cite{Gr} (using \cite[0.7]{mos}) in the case $X\eq\C^r$, $Z\eq\{0\}$, and $M\eq\OO_X(*D)$ with $D\sst X$ a hypersurface defined by a homogeneous polynomial having an isolated singularity at the origin, see \hl{2.2}{2.2} below.
\sk
In Section 1 we review the $V$-filtration of Kashiwara and Malgrange together with its relation to the Verdier specialization. In Section 2 we prove the main theorem using the assertions in the previous section and give an example.
\msn
{\bf Acknowledgement.} BD was partially supported by NSF grant DMS-1840234. MS was partially supported by JSPS Kakenhi 15K04816.
\bs\bs
\vbox{\centerline{\bf 1. Verdier specialization and $V$-filtration}
\bsn
In this section we review the $V$-filtration of Kashiwara and Malgrange together with its relation to the Verdier specialization.}
\par\htt{1.1}{}\bsn
{\bf 1.1.~$V$-filtration.} Let $Z$ be a smooth closed subvariety of a smooth complex algebraic variety $X$ with codimension $r$. Let $\I_Z$ be the ideal sheaf of $Z$ in $X$. The filtration $V$ on $\D_X$ along $Z$ is indexed by $\Z$, and is defined by
$$V^i\D_X\defs\{P\ins\D_X:P\I_Z^k\sst\I_Z^{k+i}\q\h{for any}\,\, k\ge 0\},$$
with $\I_Z^k\defs\OO_X$ for $k\less 0$.
\sk
Let $\te$ be a locally defined vector field such that $\te\ins V^{0}\D_X$ and its action on $\I_Z/\I_Z^{2}$ is the identity. Let $M$ be a quasi-unipotent regular holonomic left $\D_X$-module $M$. The $V$-filtration of Kashiwara \cite{Ka} and Malgrange \cite{Ma} of $M$ along $Z$ indexed by $\Q$ (instead of $\Z$) is an exhaustive decreasing filtration $V$ which is indexed discretely and left-continuously by $\Q$, and satisfies the following conditions (see \cite{BMS}):
\par\vbox{\htt{i}{}\msn
(i)\, the $V^{\al}M$ are coherent $V^{0}\D_X$-submodules of $M$,
\par\htt{ii}{}\skn
(ii)\, $V^i\D_XV^{\al}M\subset V^{\al+i}M$ for any $i\ins\Z,\al\ins\Q$,
\par\htt{iii}{}\skn
(iii)\, $V^i\D_XV^{\al}M\eq V^{\al+i}M$ for any $i > 0$ if $\al\gg 0$,
\par\htt{iv}{}\skn
(iv)\, the action of $\te\mi\al\pl r$ is nilpotent on $\Gr_V^{\al}M$ for any $\al\ins\Q$.}
\par\htt{1.2}{}\msn
{\bf 1.2.~ Verdier specialization.}
With the notation of \hl{1.1}{1.1}, set
$$\Xt\defs{\rm Spec}_X\bl(\mopl_{i\in\Z}\I_Z^{-i}{\otimes}t^i\br),$$
with $\I_Z^{-i}\eq\OO_X$ for $i\gess 0$. There is a natural morphism $p\co\Xt\tos\Ab_{\C}^1\defs{\rm Spec}\,\C[t]$. Its fiber over $0$ is the normal bundle
$$T_ZX\eq{\rm Spec}_X\bl(\mopl_{i\less 0}\I_Z^{-i}/\I_Z^{-i+1}\otimes t^i\br),$$
and $\Xt^*\defs p^{-1}(\Ab_{\C}^1\stm\{0\})$ is isomorphic to $X{\times}(\Ab_{\C}^1\stm\{0\})$.
So $p$ gives a deformation of $T_ZX$ to $X$, see \cite{Ve}.
\sk
Let $M$ be a quasi-unipotent regular holonomic left $\D_X$-module. Set
$$\Mt\defs\mopl_{i\in\Z}\,M{\otimes}t^i.$$
This has a structure of $\D_X\otimes_{\C}\C[t,t^{-1}]\langle\dd_t\rangle$-module, and is identified with the pull-back of $M$ by the projection $q\co\Xt^*\tos X$. We have a canonical isomorphism as $\D_X\otimes_{\C}\C[t]\langle\dd_t\rangle$-modules
\htt{1.2.1}{}
$$\rho_*j_*q^*M=\Mt.
\leqno(1.2.1)$$
Here $j\co\Xt^*\tos\Xt$, $\rho\co\Xt\tos X$ are natural affine morphisms, and $\rho_*$, $j_*$ are direct images as Zariski sheaves. Note that the direct image $j_*$ as $\D$-modules is given by the sheaf-theoretic direct image, since $j$ is an open immersion.
\sk
Let $V$ be the $V$-filtration of Kashiwara \cite{Ka} and Malgrange \cite{Ma} indexed by $\Q$ (instead of $\Z$) on $M$ along $Z$ and on $j_*q^*M$ along $T_ZX\eq p^{-1}(0)$. The {\it Verdier specialization\1} $\Sp_ZM$ of $M$ along $Z$ is defined by
\htt{1.2.2}{}
$$\Sp_ZM:=\mopl_{\al\in(0,1]}\Gr_V^{\al}(j_*q^*M)\,\bl(=\psi_t(j_*q^*M)\br).
\leqno(1.2.2)$$
By (\hl{1.2.1}{1.2.1}) its sheaf-theoretic direct image by the {\it affine\1} morphism $\rho$ is expressed as
\htt{1.2.3}{}
$$\rho_*\Sp_ZM=\mopl_{\al\in(0,1]}\Gr_V^{\al}\Mt=\mopl_{\al\in(r-1,r]}\mopl_{i\in\Z}\,\Gr_V^{\al-i}M{\otimes}t^i,
\leqno(1.2.3)$$
where the shift of $V$-filtrations comes from condition~(\hl{iv}{iv}) in \hl{1.1}{1.1} and Remark~\hl{R1.2}{1.2} below. Indeed, there are canonical isomorphisms
\htt{1.2.4}{}
$$\rho_*V^{\al}j_*q^*M=V^{\al}\Mt=\mopl_{i\in\Z}\,V^{\al+r-1-i}M{\otimes}t^i\q(\forall\,\al\ins\Q),
\leqno(1.2.4)$$
hence
\htt{1.2.5}{}
$$\rho_*\Gr_V^{\al}j_*q^*M=\Gr_V^{\al}\Mt =\mopl_{i\in\Z}\,\Gr_V^{\al+r-1-i}M{\otimes}t^i\q(\forall\,\al\ins\Q),
\leqno(1.2.5)$$
see \cite[1.3]{BMS}. These are crucial to the relation between the construction of Kashiwara \cite{Ka} and the Verdier specialization \cite{Ve}.
\sk
Assume now $M$ underlies a mixed Hodge module. It has the Hodge filtration $F$, and
\htt{1.2.6}{}
$$\rho_*j_*q^*F_pM=\mopl_{i\in\Z}\,F_pM{\otimes}t^i\q(\forall\,\al\ins\Q),
\leqno(1.2.6)$$
since the composition $\rho\ssc j$ coincides with the projection $X{\times}(\Ab_{\C}^1\stm\{0\})\tos X$. (Note that the Hodge filtration on the underlying {\it left\1} $\D$-module is {\it not\1} shifted by the smooth pullback $q^*$; consider for instance the structure sheaf case.)
\sk
By \cite[Prop.\,3.2.2]{mhp} together with (\hl{1.2.4}{1.2.4}) we then get that
\htt{1.2.7}{}
$$\aligned\rho_*F_pV^{\al}(j_*q^*M)&=\rho_*j_*q^*F_pM\cap\rho_*V^{\al}j_*q^*M\\&=\mopl_{i\in\Z}\,F_pV^{\al+r-1-i}M{\otimes}t^i\q\q(\forall\,\al\sgt0),\endaligned
\leqno(1.2.7)$$
since $\rho$ is an affine morphism. Hence
\htt{1.2.8}{}
$$\rho_*F_p\Gr_V^{\be}(j_*q^*M)=\mopl_{i\in\Z}\,F_p\Gr_V^{\be+r-1-i}M{\otimes}t^i\q\q(\forall\,\be\ins(0,1]).
\leqno(1.2.8)$$
This implies (\hl{6}{6}) in the introduction using Remark~\hl{R1.2}{1.2} below, see also \cite[Lem.\,2.4]{CD}.
\par\htt{R1.2}{}\msn
{\bf Remark~1.2.} Let $x_1,\dots,x_n$ be local coordinates of $X$ such that $Z\eq\mcap_{i=1}^r\{x_i\eq0\}$ locally. Then $\Xt$ has local coordinates $\xt_1,\dots,\xt_n,\tti$ such that $\xt_i\eq x_i/t\,\,(i\less r)$, $\xt_i\eq x_i\,\,(i\sgt r)$, $\tti\eq t$ on $\Xt^*$, and we have the equalities
\htt{1.2.9}{}
$$t\dd_t =\tti\dd_{\tti}-\msum_{i=1}^r\,\xt_i\dd_{\xt_i},\q x_i\dd_{x_i}=\xt_i\dd_{\xt_i}\,\,\,(\forall\,i).
\leqno(1.2.9)$$
Since $M{\otimes}t^i$ is identified with ${\rm Ker}(t\dd_t\mi i)\sst\Mt$, we see that
\htt{1.2.10}{}
$$\h{$\tti\dd_{\tti}\mi i$ is identified with $\te\defs\msum_{i=1}^r\,\xt_i\dd_{\xt_i}$ on $M{\otimes} t^i$.}
\leqno(1.2.10)$$
The existence of the filtration $V$ can be reduced to the hypersurface case by this argument, see \cite{BMS}.
\sk
The following is noted in \cite[2.30]{mhm} in case $Z\sst X$ is defined by {\it global\1} functions on $X$.
\par\htt{P1.2}{}\msn
{\bf Proposition~1.2.} {\it We have the isomorphisms
\htt{1.2.11}{}
$$i^*_Z\M=(i'_Z)^*\Sp_Z\M,\q i^!_Z\M=(i'_Z)^!\Sp_Z\M,
\leqno(1.2.11)$$
where $i_Z\co Z\into X$ and $i'_Z\co Z\into T_ZX$ are natural inclusions.}
\msn
{\it Proof.} This can be reduced to
\htt{1.2.12}{}
$$\aligned\Sp_Z(j_{X\setminus Z})_!j^*_{X\setminus Z}\M&=(j_{X'\setminus Z})_!j^*_{X'\setminus Z}\Sp_Z\M,\\ \Sp_Z(j_{X\setminus Z})_*j^*_{X\setminus Z}\M&=(j_{X'\setminus Z})_*j^*_{X'\setminus Z}\Sp_Z\M,\endaligned
\leqno(1.2.12)$$
where $j_{X\setminus Z}\co X\stm Z\into X$ and $j_{X'\setminus Z}\co X'\stm Z\into X'$ are natural inclusions with $X'\defs T_ZX$. Indeed, we have the distinguished triangles
\htt{1.2.13}{}
$$\aligned(j_{X\setminus Z})_!j^*_{X\setminus Z}\to{\rm id}\to(i_Z)_*i^*_Z\buildrel{[1]}\over\to,\\(i_Z)_*i^!_Z\to{\rm id}\to(j_{X\setminus Z})_*j^*_{X\setminus Z}\buildrel{[1]}\over\to,\endaligned
\leqno(1.2.13)$$
and similarly with $X$, $i_Z$ replaced by $X'$, $i'_Z$ respectively. Moreover $\Sp_Z$ is the identity on mixed Hodge modules supported on $Z\sst X$.
\sk
The isomorphisms in (\hl{1.2.12}{1.2.12}) are reduced to the corresponding ones for the underlying $\Q$-complexes using the functorial morphisms ${\rm id}\to(j_{X\setminus Z})_*j_{X\setminus Z}^*$, etc.\ together with the mapping cone. (Note that a bounded complex of mixed Hodge modules is acyclic if its underlying $\Q$-complex is.) The desired isomorphisms for the underlying $\Q$-complexes can be shown easily by the definition of Verdier specialization in the case of open direct images with proper supports. We can then apply the commutativity of the Verdier specialization with the dual. This finishes the proof of Proposition~\hl{P1.2}{1.2}.
\bs\bs
\vbox{\centerline{\bf 2. Proof of the main theorem}
\bsn
In this section we prove the main theorem using the assertions in the previous section and give an example.}
\par\htt{2.1}{}\msn
{\bf 2.1.~Proof of Theorem~\hl{T1}{1}.} We prove the assertion together with (\hl{8}{8})--(\hl{9}{9}) by increasing induction on $r\gess1$. We may assume $r\gess2$, since they hold for $r\eq1$ by the theory of Hodge modules, see \cite{mhp}, \cite{ypg}. We also assume that the vector bundle is trivialized by $t_1,\dots,t_r$ as in the introduction, since the assertion is local. Here we may assume that $t_r$ is sufficiently general replacing it with a sufficiently general $\C$-linear combination of $t_1,\dots,t_r$ (since only the trivialization of vector bundle is changed).
\sk
Let $V_r$ be the $V$-filtration on $\D_X$, $M'$ along $\{t_r\eq0\}$ indexed by $\Z$ and $\Q$ respectively. Since $\te\ins V_r^0\D_X$, the $V_r^{\al}M'$ are stable by $\te$ for any $\al\ins\Q$. The filtration $V_r$ is then compatible with the direct sum decomposition (\hl{2}{2}) using \cite[Rem.\,A.7c]{rh}, since any local section of $M'$ is contained in a finite direct sum in (\hl{2}{2}), and is annihilated by a polynomial of $\theta$. This gives a filtration (also denoted by $V_r$) on $K_{\!Z}^{!,\ssb}(M')_{\al}$, $K_{\!Z}^{*,\ssb}(M')_{\al}$ so that there are isomorphisms of complexes of filtered $\D_Z$-modules
\htt{2.1.1}{}
$$\aligned&\Gr_{V_r}^{\be}K_{\!Z}^{!,\ssb}(M',F)_{\al}\\&=\bl[K_{\!Z}^{!,\ssb}(\Gr_{V_r}^{\be}(M',F))_{\al-\be}\buildrel{t_r\,}\over\to K_{\!Z}^{!,\ssb}(\Gr_{V_r}^{\be+1}(M',F))_{\al-\be}\br],\\&\Gr_{V_r}^{\be}K_{\!Z}^{*,\ssb}(M',F)_{\al}\\&=\bl[K_{\!Z}^{*,\ssb}(\Gr_{V_r}^{\be+1}(M',F))_{\al-\be}\buildrel{\dd_r\,}\over\to K_{\!Z}^{*,\ssb}(\Gr_{V_r}^{\be}(M',F[-1]))_{\al-\be}\br].\endaligned
\leqno(2.1.1)$$
Here we denote by $[A\tos B]$ the shifted mapping cone $C(A\tos B)[-1]$, and $K_{\!Z}^{!,\ssb}(\Gr_{V_r}^{\be}M')$ and $K_{\!Z}^{*,\ssb}(\Gr_{V_r}^{\be}M')$ are respectively the Koszul complexes for the actions of $t_1,\dots,t_{r-1}$ and $\dd_1,\dots,\dd_{r-1}$ on the $\D_{X'}$-module $\Gr_{V_r}^{\be}M'$ with $X'\defs\{t_r\eq0\}\sst X$, where the grading is defined by using $\te'\defs\msum_{i=1}^{r-1}\,t_i\dd_i$ (hence $\te\eq\te'{+}t_r\dd_r$). Note that the Koszul complex for the actions of $t_1,\dots,t_r$ is identified with the shifted mapping cone of the action of $t_r$ on the Koszul complex for $t_1,\dots,t_{r-1}$, and similarly with $t_1,\dots,t_r$ replaced by $\dd_1,\dots,\dd_r$.
\sk
By the relation to the Hodge filtration \cite[3.2.1]{mhp}, we have the filtered isomorphisms
\htt{2.1.2}{}
$$\aligned t_r:(V_r^{\be}M',F)\simto(V_r^{\be+1}M',F)&\q\q\h{if}\q\be\sgt0,\\ t_r:\Gr_{V_r}^{\be}(M',F)\simto\Gr_{V_r}^{\be+1}(M',F)&\q\q\h{if}\q\be\sgt0,\\ \dd_r:\Gr_{V_r}^{\be+1}(M',F)\simto\Gr_{V_r}^{\be}(M',F[-1])&\q\q\h{if}\q\be\slt0,\endaligned
\leqno(2.1.2)$$
as well as the isomorphisms (using condition~(\hl{iv}{iv}) in \hl{1.1}{1.1})
\htt{2.1.3}{}
$$\aligned t_r:\Gr_{V_r}^{\be}M'\simto\Gr_{V_r}^{\be+1}M'&\q\q\h{if}\q\be\slt0,\\
\dd_r:\Gr_{V_r}^{\be+1}M'\simto\Gr_{V_r}^{\be}M'&\q\q\h{if}\q\be\sgt0.\endaligned
\leqno(2.1.3)$$
Using (\hl{2.1.1}{2.1.1}) together with (\hl{2.1.2}{2.1.2}) and the inductive hypothesis, we can verify that the filtered complexes
$$\Gr_{V_r}^{\be}K_{\!Z}^{!,\ssb}(M',F)_{\al},\q\Gr_{V_r}^{\be}K_{\!Z}^{*,\ssb}(M',F)_{\al},\q V_r^{\be}K_{\!Z}^{!,\ssb}(M',F)_{\al}\q V_r^{\be}K_{\!Z}^{*,\ssb}(M',F)_{\al}$$
are filtered acyclic if the following conditions are satisfied respectively:
$$\be\sgt0\,\,\h{or}\,\,\al\sgt\be,\q\be\slt0\,\,\h{or}\,\,\al\slt\be,\q\be\sgt0,\q\al\slt\be.$$
Indeed, for $\Gr_{V_r}^{\be}K_{\!Z}^{!,\ssb}(M',F)_{\al}$, we use (\hl{2.1.1}{2.1.1}) together with (\hl{2.1.2}{2.1.2}) if $\be\sgt0$, and the inductive hypothesis if $\al{-}\be\sgt0$. (If $\be\notins\Z$, we consider a mixed Hodge module which is isomorphic to the non-unipotent monodromy part of the nearby cycle Hodge module $\psi_{t_r,\ne1}\M'$ up to a Tate twist and whose underlying filtered $\D$-module is the direct sum of $\Gr_{V_r}^{\be'}(M',F)$ for $\be'\ins\Q\cap(\lfloor\be\rfloor,\lceil\be\rceil)$ in order to apply the inductive hypothesis.) For $\Gr_{V_r}^{\be}K_{\!Z}^{*,\ssb}(M',F)_{\al}$, the argument is similar with inequalities reversed. For $V_r^{\be}K_{\!Z}^{!,\ssb}(M',F)_{\al}$, we consider the mapping cone similar to the first one in (\hl{2.1.1}{2.1.1}) with $\Gr_{Vr}^{\be}$ replaced by $V_r^{\be}$ (where the grading by the action of $\te$ instead of $\te'$ is used), and apply the first isomorphism of (\hl{2.1.2}{2.1.2}). Note that the latter isomorphism is compatible with the direct sum decomposition by the eigenvalues of the action of $\te$. For $V_r^{\be}K_{\!Z}^{*,\ssb}(M',F)_{\al}$, we apply Lemma~\hl{L2.1}{2.1} below.
\sk
We then see that (\hl{8}{8})--(\hl{9}{9}) hold for $r$. We get also the following isomorphisms in the bounded derived category of $F$-filtered $\D_Z$-modules $D^bF(\D_Z)$
\htt{2.1.4}{}
$$\aligned K_{\!Z}^{!,\ssb}(M')_0&=\Gr_{V_r}^0K_{\!Z}^{!,\ssb}(M')_0\\&=\bl[K_{\!Z}^{!,\ssb}(\Gr_{V_r}^0M')_0\buildrel{t_r\,}\over\to K_{\!Z}^{!,\ssb}(\Gr_{V_r}^1M')_0\br],\\ K_{\!Z}^{*,\ssb}(M')_0&=\Gr_{V_r}^0K_{\!Z}^{*,\ssb}(M')_0\\&=\bl[K_{\!Z}^{*,\ssb}(\Gr_{V_r}^1M')_0\buildrel{\dd_r\,}\over\to K_{\!Z}^{*,\ssb}(\Gr_{V_r}^0M')_0\br].\endaligned
\leqno(2.1.4)$$
Here $\Gr_{V_r}^0M'$ is supported on $Z$ (shrinking $Z$ if necessary) since $t_r$ is sufficiently general, see Remark~\hl{R2.1b}{2.1b} below. We then get that
\htt{2.1.5}{}
$$\aligned\Hc^jK_{\!Z}^{!,\ssb}(\Gr_{V_r}^0M')&\eq0\q\h{if}\,\,\,j\nes0,\\ \Hc^jK_{\!Z}^{*,\ssb}(\Gr_{V_r}^0M')&\eq0\q\h{if}\,\,\,j\nes r{-}1.\endaligned
\leqno(2.1.5)$$
Here $t_r$ and $\dd_r$ in the mapping cones underlie the morphisms can and Var between the nearby and vanishing cycle functors of mixed Hodge modules, which are strictly compatible with the Hodge filtration. So the assertion follows from the inductive hypothesis using \cite[Remark after Cor.\,2.24]{mhm} and Remark~\hl{R2.1c}{2.1c} below. The assertion on the weight filtration follows from Proposition~\hl{P2.1}{2.1} and Remark~\hl{R2.1a}{2.1a} below. This finishes the proof of Theorem~\hl{T1}{1}.
\par\htt{R2.1a}{}\msn
{\bf Remark~2.1a.} Let $M$ be the underlying $\D_X$-module of a monodromical {\it pure\1} Hodge module $\M'$ on a vector bundle $X\tos Z$. Then the action of $\te\mi\al\pl r$ on $M'_{\al}$ in (\hl{1}{1}) is {\it semisimple.} Indeed, a monodromical Hodge module is invariant by the Verdier specialization. So the assertion follows from (\hl{1.2.10}{1.2.10}), since the weight filtration on $\psi_{\tti,\la}$ is given by the monodromy filtration for the action of $\tti\dd_{\1\tti}\mi\al$ for $\al\in(-1,0]$ with $\la\eq e^{-2\pi i\al}$.
\par\htt{R2.1b}{}\msn
{\bf Remark~2.1b.} Let $\M'$ be a monodromical mixed Hodge module on a trivial vector bundle $\pi\co X\tos Z$, which is trivialized by sections $t_1,\dots,t_r$ as in the introduction. Take $p\ins Z$. Assume $t_r$ is sufficiently general replacing it with a sufficiently general $\C$-linear combination of $t_1,\dots,t_r$. Then the vanishing cycle Hodge module $\varphi_{t_r}\M'$ is supported on $Z$ replacing $X$ by $\pi^{-1}(U)$ with $U$ a sufficiently small neighborhood of $p$ in $Z$. This can be verified by taking a Whitney stratification of the trivial projective bundle $P_Z:=(X\stm Z)/\C^*$ which is compatible with the monodromical Hodge module $\M'$ and also with the fiber $P_{Z,p}$ at $p$. We choose $t_r$ so that $\{t_r\eq0\}\sst P_{Z,p}$ is transversal to any stratum contained in $P_{Z,p}$, and use the Whitney's condition (a) on the limits of tangent spaces.
\par\htt{R2.1c}{}\msn
{\bf Remark~2.1c.} Let $\phi\co(A^{\ssb},F)\tos(B^{\ssb},F)$ be a morphism of filtered complexes in an abelian category $\A$. Let $a$ be a positive integer. Assume $A^j\eq B^j\eq0$ for $j\,{\notin}\,\{0,\dots,a\}$, $(A^{\ssb},F)$, $(B^{\ssb},F)$ are strict filtered complexes, the filtered morphisms $H^j\phi\co H^j(A^{\ssb},F)\tos H^j(B^{\ssb},F)$ are strict for $j\eq0,a$, and either $H^jA^{\ssb}\eq0$ for any $j\nes0$ or $H^jB^{\ssb}\eq0$ for any $j\nes a$. Then the mapping cone of $\phi\co(A^{\ssb},F)\tos(B^{\ssb},F)$ is a strict filtered complex. (This can be verified for instance using an abelian category whose objects are direct sums of objects of $\A$, and considering the functor associating $\mopl_p\,F_pA$ to an filtered object $(A,F)$, see also \cite[1.3]{mhp}.)
\par\htt{R2.1d}{}\msn
{\bf Remark 2.1d.} Let $t'_1,\dots,t'_r$ be local sections of $X\tos Z$ given by another local trivialization of vector bundle so that $t'_i\eq\msum_{j=1}^r\,a_{i,j}t_j$ with $a_{i,j}$ local sections of $\OO_Z$. Let $z_1,\dots,z_n$ be local coordinates of $Z$. Setting $z'_k\eq z_k$ ($k\ins[1,n]$), we get two local coordinate systems $z_k,t_j$ and $z'_k,t'_j$ of $X$, and
\htt{2.1.6}{}
$$\aligned\dd_{z_k}&\eq\dd_{z'_k}+\msum_{i=1}^r\msum_{j=1}^r\,(\dd_{z_k}a_{i,j})\1t_j\dd_{t'_i}\q(k\ins[1,n]),\\ \dd_{t_j}&\eq\msum_{i=1}^r\,a_{i,j}\1\dd_{t'_i}\q(j\ins[1,r]).\endaligned
\leqno(2.1.6)$$
So the ambiguities of liftings of vector fields on $Z$ to $X$ are given by $\OO_Z$-linear combinations of the $t_i\dd_j$ when a local trivialization is changed.
\par\htt{R2.1e}{}\msn
{\bf Remark 2.1e.} In the notation of Theorem~\hl{T1}{1}, let $\Lc^{\ssb}\tos M'$ be a flat resolution as $\D_X$-module. (There is a locally free resolution globally only in the quasi-projective variety case.) Using the tensor product with the Koszul complex associated with the left multiplications by $t_k$ ($k\ins[1,r]$) on $\D_X$, we get the vanishing
\htt{2.1.7}{}
$$\Hc^jK_{\!Z}^{!,\ssb}(\Lc^i)\eq0\q\h{if}\,\,\,j\nes r,
\leqno(2.1.7)$$
since the tensor product with a flat module is an {\it exact functor.} Moreover the action of a vector field on $Z$ is well defined on $\Hc^rK_{\!Z}^{!,\ssb}(\Lc^i)$ by Remark~\hl{R2.1d}{2.1d}. So we get a {\it globally well defined\1} complex of $\D_Z$-modules $K_{\!Z}^{!,\ssb}(M')^{\bf L}$ whose $i\1$th component is $\Hc^rK_{\!Z}^{!,\ssb}(\Lc^{i-r})$. (Recall that $\Lc^i\eq0$ unless $i\less0$.) This is compatible with the definition $\D_{Z\to X}\defs\OO_Z{\otimes}_{\OO_X}\D_X$.
\sk
Using a standard argument about double complexes (or a spectral sequence related to it), we get a canonical isomorphism of $\D_U$-modules
\htt{2.1.8}{}
$$\Hc^jK_{\!Z}^{!,\ssb}(M')^{\bf L}|_U\eq\Hc^jK_{\!Z}^{!,\ssb}(M')|_U,
\leqno(2.1.8)$$
where $U\sst Z$ is an open subset on which local sections $t_j$ are defined (and the action of vector fields on the right-hand side is defined by using these $t_j$). This says that the $\D_Z$-module structure on each $\Hc^jK_{\!Z}^{!,\ssb}(M')$ is defined independently of the choice of local sections $t_j$. Here we can use the flatness of $\D_X$ over $\OO_X$ to prove that flatness over $\D_X$ implies flatness over $\OO_X$, and show that the isomorphism is compatible with the canonical isomorphism for the underlying $\OO$-modules (which are defined without choosing the $t_j$) by using the tensor product with the Koszul complex for the left multiplications by $t_k$ ($k\ins[1,r]$) on $\OO_X$.
\par\htt{R2.1f}{}\msn
{\bf Remark 2.1f.} For $K_{\!Z}^{*,\ssb}(M')$, the argument is similar to Remark~\hl{R2.1e}{2.1e}. We replace a flat resolution $\Lc^{\ssb}\tos M'$ by a (non-quasi-coherent) injective resolution $M'\tos\I^{\ssb}$ (and ${\bf L}$ by ${\bf R}$). We have $\Hc^jK_{\!Z}^{*,\ssb}(\I^i)\eq0$ if $j\nes0$. Indeed, the restriction of $\I^i$ to any open subset is injective (using the zero-extension, which is not quasi-coherent), and the functor assigning $\Hom_{\pi_*\D_X}(\Nc,\pi_*\I^i)$ to any {\it quasi-coherent\1} $\pi_*\D_X$-module $\Nc$ gives an exact functor, since $\pi$ is an affine morphism and there is a quasi-coherent $\D_X$-module $\Nc'$ with $\pi_*\Nc'\eq\Nc$, where $\D_X{\otimes}_{\pi^{-1}\pi_*\D_X}\pi^{-1}\Nc\eq\Nc'$ using quasi-coherence (via local presentations) and the {\it right exactness\1} of tensor product. We consider the Koszul complex for the right multiplications by $\dd_k$ ($k\ins[1,r]$) on $\pi_*\D_X$, and also on the polynomial ring of the $\dd_k$ with coefficients in $\OO_Z$ (over which $\pi_*\D_X$ is flat). The latter subring corresponds to $\OO_X\sst\D_X$ in the flat resolution case.
\sk
For the compatibility with the definition of direct image for $\D$-modules, note that $\D_{X\to Z}$ is locally isomorphic to $\D_X/\msum_{k=1}^r\,\D_X\dd_k$, and it is twisted by $\om_{X/Z}$ to get $\D_{Z\leftarrow X}$. These produce the relative de Rham complex as in (\hl{4}{4})--(\hl{5}{5}) using the Koszul complex for the actions of $\dd_k$ on $\pi_*M'$. Here it is enough to consider its initial cohomology sheaf, that is, the intersection of the kernels of $\dd_k$, if we replace $M'$ with $\I^i$.
\par\htt{R2.1g}{}\msn
{\bf Remark~2.1g.} Set $R\defs V_r^0\Gr^F_{\ssb}\Gr_V^0\D_{X,p}$, $S\defs\Gr^F_{\ssb}\D_{Z,p}$, where $p\ins Z\sst X$. These are graded $\OO_{Z,p}$-algebras, and $R$ is generated over $S$ by
$$\xi_{i,j}\defs\Gr^F_1t_i\dd_j\ins R_1\q\h{for}\,\,\,(i,j)\ins[1,r]^2\stm[1,r{-}1]{\times}\{r\}.$$
The $V_r^{\be}\Gr^F_{\ssb}M'_{\al,p}$ are finite graded $R$-modules. Set $\eta_j\defs\Gr^F_1\dd_j$.
\sk
Let $\be$ be a positive rational number. We can verify that $V_r^{\be}\Gr^F_{\ssb}M'_p$ is a finite graded module over $V_r^0\Gr^F_{\ssb}\D_{X,p}$ by an argument similar to the proof of \cite[Lem.\,1.1]{JKSY} (using Nakayama's lemma). We then have a sufficiently small rational number $\al_0$ such that
\htt{2.1.9}{}
$$V_r^{\be}\Gr^F_{\ssb}M'_{\al,p}=\msum_{j=1}^{r-1}\,\eta_j\1V_r^{\be}\Gr^F_{\ssb\1-1}M'_{\al+1,p}\q\h{if}\,\,\,\be\sgt0,\,\,\al{-}\be\less\al_0,
\leqno(2.1.9)$$
reducing to the case $\be\ins(0,1]$ by using the first isomorphism of (\hl{2.1.2}{2.1.2}).
\sk
By (\hl{2.1.2}{2.1.2}) and (\hl{2.1.9}{2.1.9}) we get that
\htt{2.1.10}{}
$$V_r^{\be}\Gr^F_{\ssb}M'_{\al,p}=\msum_{j=1}^{r-1}\,\xi_{r,j}\1V_r^{\be-1}\Gr^F_{\ssb\1-1}M'_{\al,p}\q\h{if}\,\,\,\be\sgt1,\,\,\al{-}\be\less\al_0.
\leqno(2.1.10)$$
Note that the last condition may hold if $\be$ is sufficiently large for {\it any\1} $\al\ins\Q$. The equality (\hl{2.1.10}{2.1.10}) says that the filtration $V_r$ on $\Gr^F_{\ssb}M'_{\al,p}$ is essentially the $I$-adic filtration, where $I\sst R$ is the graded ideal generated by the $\xi_{r,j}$ for $j\ins[1,r{-}1]$.
\par\htt{L2.1}{}\msn
{\bf Lemma~2.1.} {\it If $\Gr_{V_r}^{\be}K_{\!Z}^{*,\ssb}(M',F)_{\al}$ is filtered acyclic for all $\be\gess\be_0$, then $V_r^{\be_0}K_{\!Z}^{*,\ssb}(M',F)_{\al}$ is filtered acyclic $($with $\al\ins\Q$ fixed$\1)$.}
\msn
{\it Proof.} Using the Mittag-Leffler condition (see \cite[Prop.\,13.2.3]{Gro}), the completion of the complex $V_r^{\be_0}\Gr^F_{\ssb}K_{\!Z}^{*,\ssb}(M'_p)_{\al}$ by the filtration $V_r$ is acyclic, where $p\ins Z$. This completion is the same as the $I$-adic completion by Remark~\hl{R2.1g}{2.1g}, and can be given by the tensor product with the $I$-adic completion $\widehat{R}$ of $R$, where $\widehat{R}$ is flat over $R$, see \cite[8.7--8]{Mat}. We thus get the vanishing of the $I$-adic completion of
$$H^jV_r^{\be_0}\Gr^F_{\ssb}K_{\!Z}^{*,\ssb}(M'_p)_{\al}.$$
The last assertion holds with $I$ replaced by the maximal homogeneous ideal $I'\sst R$ with $R/I'\eq\OO_{Z,p}$ using the tensor product with the $I'$-adic completion of $R$. So the cohomology graded $R$-modules vanish without taking the completion. (Here we can use also the graded Nakayama lemma.) This finishes the proof of Lemma~\hl{L2.1}{2.1}.
\par\htt{P2.1}{}\msn
{\bf Proposition~2.1.} {\it Let $\M'$ be a monodromical pure Hodge module of weight $w$ on a vector bundle $X\tos Z$. Then the $H^ji_Z^*\M'$ and $H^ji_Z^!\M'$ are pure of weight $w{+}j$ for $j\ins\Z$.}
\msn
{\it Proof.} We may assume that $\M'$ is simple and $Y\defs{\rm Supp}\,\M'$ is not contained in $Z$, since the assertion is clear otherwise. In this proof, we denote by $\rho:\Xt\tos X$ the blow-up of $X$ along $Z$ with $i_E\co E\into\Xt$ the exceptional divisor. Note that $\Xt$ is a line bundle over $E$ with structure morphism $q\co \Xt\to E$. There is a unique simple pure Hodge module $\Mct'$ on $\Xt$ such that $\Mct'|_{\Xt\setminus E}\eq\M'|_{X\setminus Z}$. Let $\Yt\sst\Xt$ be the proper transform of $Y$. The local monodromy around $\Yt\cap E$ of the generic local system of $\Mct'$ is semisimple using a non-characteristic restriction of $\M'$ to a sufficiently general line contained in the support of $\M'$ and also in a fiber of $\pi$. (Here we first take the non-characteristic restriction to a fiber $X_p$ of $X\to Z$, and then take $d{-}1$ general hyperplanes containing the line, where $d$ is the dimension of the intersection of the support of $\M'$ with $X_p$. Note that the semisimplicity of the action of the Euler vector field $\te$ on an algebraic regular holonomic $\D_{\C^*}$-module $\M''$ which is finite over $\OO_{\C^*}$ is equivalent to the semisimplicity of the corresponding local system on $\C^*$. Here the first semisimplicity means that $\Gamma(\C^*,\M'')$ is a direct sum of the kernels of $\te\mi\al$ for $\al\in\C$.)
\sk
Using the compatibility of intermediate direct images \cite{BBD} with external products, we see that the monodromy of the nearby cycle Hodge module $\psi_g\Mct'$ along a local defining function $g$ of $E\sst\Xt$ is also semisimple. This implies that the unipotent monodromy part of the vanishing cycle Hodge module $\varphi_{g,1}\Mct'$ vanishes, hence $H^{-1}i_E^*\Mct'$ and $H^1i_E^!\Mct'$ are pure Hodge modules of weight $w-1$ and $w+1$ respectively, and moreover $H^ki_E^*\Mct'\eq0$ for $k\nes{-}1$ and $H^ki_E^!\Mct'\nes0$ for $k\nes1$. Since $\M'$ is a direct factor of the direct image $\rho_*\Mct'$ (using for instance \cite[(4.5.2), (4.5.4)]{mhm}), the assertion follows from the base change theorem and the stability of pure complexes under the direct images by proper morphisms, see for instance \cite[(4.4.3), (4.5.2)]{mhm}. This finishes the proof of Proposition~\hl{P2.1}{2.1}.
\par\htt{2.2}{}\msn
{\bf 2.2.~Example.} Assume $X\eq\C^r$, $Z\eq\{0\}$, and $M'\eq\OO_X(*D)$, where $D\sst X$ is a hypersurface defined by a homogeneous polynomial $f$ having an isolated singularity at 0 with $r\gess 3$. Set $U\defs\PP^{r-1}\stm\{f\eq0\}$. Then $X\stm D$ is a $\C^*$-bundle over $U$, and we have the isomorphisms of mixed Hodge structures
\htt{2.2.1}{}
$$H^j(X\stm D,\Q)\cong\begin{cases}H^{r-1}(U,\Q)(-1)&\h{if}\,\,\,\,j\eq r,\\H^{r-1}(U,\Q)&\h{if}\,\,\,\,j\eq r{-}1,\\ \Q(-1)&\h{if}\,\,\,\,j\eq 1,\\ \Q&\h{if}\,\,\,\,j\eq 0,\\ \,0&\h{otherwise},\end{cases}
\leqno(2.2.1)$$
see for instance \cite[1.3]{RSW}. (Here the Thom class vanishes with $\Q$-coefficients.)
\sk
We see that $\Hc^jK_{\!Z}^{*,\ssb}(M')_0$ is generated by $1\ins\OO_{X,0}$ if $j\eq0$, and $\ddd f/f\ins\Om_{X,0}(*D)$ if $j\eq1$. Here the pullback $i_0^*$ is identified with the direct image $\pi_*$ of a monodromical $\D$-module, and $\D_Z$-modules are identified with $\C$-vector spaces. Set
$$\Om_f\defs\Om_{X,0}^r/\ddd f{\wedge}\Om_{X,0}^{r-1}.$$
This is a finite-dimensional $\Z$-graded $\C$-vector space, and is well studied in the hypersurface isolated singularity theory. By \cite{Gr} we have the isomorphisms
\htt{2.2.2}{}
$$\Gr_F^{r-1-p}H^{r-1}(U,\C)=f^{-p-1}\1\Om_{f,(p+1)d}\q\q(\forall\,p\in\Z).
\leqno(2.2.2)$$
For $j\eq r$, it follows from Theorem~\hl{T1}{1} together with (\hl{2.2.1}{2.2.1}) and \cite[0.7]{mos} that
\htt{2.2.3}{}
$$\Gr^F_{p-r}\Hc^rK_{\!Z}^{*,\ssb}(M')_0\cong f^{-p-1}\1\Om_{f,(p+1)d}\q\q(\forall\,p\ins\Z),
\leqno(2.2.3)$$
where $F_p\eq F^{-p}$. For $j\eq r{-}1$, it is expected that
\htt{2.2.4}{}
$$\Gr^F_{p-r+1}\Hc^{r-1}K_{\!Z}^{*,\ssb}(M')_0\cong f^{-p-1}\1\iota_{\te}\,\1\Om_{f,(p+1)d}\q(\forall\,p\ins\Z),
\leqno(2.2.4)$$
where $\iota_{\te}$ denotes the interior product with the Euler field $\te$, see \cite[1.4]{BS}. The latter is defined by taking representatives of elements of $\Om_{f,(p+1)d}$, and similarly for $f^{-p-1}$.

\sk
{\smaller\smaller 
Department of Mathematics, University of Michigan, 530 Church Street, Ann Arbor, MI 48109, USA

{\it Email address}\,: qyc@umich.edu

Department of Mathematics, University of Michigan, 530 Church Street, Ann Arbor, MI 48109, USA

{\it Email address}\,: bdirks@umich.edu

RIMS Kyoto University, Kyoto 606-8502 Japan

{\it Email address}\,: msaito@kurims.kyoto-u.ac.jp}

\begin{thebibliography}{RSW\,21}
\bibitem[BaSa\,22]{BS} Bath, D., Saito, M., Twisted logarithmic complexes of positively weighted homogeneous divisors (arxiv:2203.11716).
\bibitem[BBD\,82]{BBD} Beilinson, A., Bernstein, J., Deligne, P., Faisceaux pervers, Ast\'erisque 100, Soc.\ Math.\ France, Paris, 1982.
\bibitem[BMS\,06]{BMS} Budur, N., Musta\c{t}\u{a}, M., Saito, M., Bernstein-Sato polynomials of arbitrary varieties, Compos.\ Math.\ 142 (2006), 779--797.
\bibitem[ChDi\,21]{CD} Chen, Q., Dirks, B., On $V$-filtration, Hodge filtration and Fourier transform (arxiv:2111.04622) to appear in Selecta Math.
\bibitem[Gri\,69]{Gr} Griffiths, P., On the period of certain rational integrals I, II, Ann.\ Math.\ 90 (1969), 460--541.
\bibitem[Gro\,61]{Gro} Grothendieck, A., El\'ements de g\'eom\'etrie alg\'ebrique III-1, Publ.\ Math.\ IHES, 11 (1961), 5--167.
\bibitem[JKSY\,22]{JKSY} Jung, S.-J., Kim, I.-K., Saito, M., Yoon, Y., Brian\c{c}on-Skoda exponents and the maximal root of reduced Bernstein-Sato polynomials, Selecta Math.\ (N.S.) 28 (2022), Paper No. 78.
\bibitem[Kas\,83]{Ka} Kashiwara, M., Vanishing cycle sheaves and holonomic systems of differential equations, Lect.\ Notes Math.\ 1016, Springer, Berlin, 1983, pp. 136--142.
\bibitem[Mal\,83]{Ma} Malgrange, B., Polyn\^ome de Bernstein-Sato et cohomologie \'evanescente, Ast\'erisque 101-102 (1983), 243--267.
\bibitem[Mat\,83]{Mat} Matsumura, H., Commutative ring theory, Cambridge University Press, 1986.
\bibitem[RSW\,21]{RSW} Reichelt, T., Saito, M., Walther, U., Dependence of Lyubeznik numbers of cones of projective schemes on projective embeddings, Selecta Math.\ (N.S.) 27 (2021), Paper No.\ 6.
\bibitem[Sa\,88]{mhp} Saito, M., Modules de Hodge polarisables, Publ.\ RIMS, Kyoto Univ.\ 24 (1988), 849--995.
\bibitem[Sa\,90]{mhm} Saito, M., Mixed Hodge modules, Publ.\ RIMS, Kyoto Univ.\ 26 (1990), 221--333.
\bibitem[Sa\,09]{mos} Saito, M., On the Hodge filtration of Hodge modules, Moscow Math.\ J.\ 9 (2009), 161--191.
\bibitem[Sa\,17]{ypg} Saito, M., A young person's guide to mixed Hodge modules, in Hodge theory and $L^2$-analysis, Adv.\ Lect.\ Math.\ 39, Int.\ Press, Somerville, MA, 2017, 517--553 (arxiv:1605.00435).
\bibitem[Sa\,22]{rh} Saito, M., Notes on regular holonomic $D$-modules for algebraic geometers (arxiv:2201.01507).
\bibitem[Ve\,83]{Ve} Verdier, J.-L., Sp\'ecialisation de faisceaux et monodromie mod\'er\'ee, Ast\'erisque, 101-102, Soc.\ Math.\ France, Paris 1983, pp.~332–364. 
\end{thebibliography}
\end{document}